\documentclass[dvipdfm]{article}
\usepackage[T1]{fontenc}
\usepackage[latin1]{inputenc}
\usepackage{amsfonts,amsmath,amssymb,amsthm}
\usepackage{graphicx}
\newtheorem{defi}{Definition}
\newtheorem{prop}{Proposition}
\newtheorem{thme}{Theorem}
\newtheorem{lemm}{Lemma}

\newcommand{\fin}{\hfill \framebox[2mm]{\ } \medskip}
\newcommand{\N}{\mathbb N}

\newcommand{\R}{\mathbb R}

\title{Asymptotic expansion of planar canard solutions near a non-generic turning point}
\author{Thomas Forget \thanks{Laboratoire de Math\'ematiques et Applications, P\^ole Sciences et Technologies - Universit\'e de La Rochelle, Avenue Michel Cr\'epeau, 17042 LA ROCHELLE - FRANCE}\\thomas.forget@univ-lr.fr}

\begin{document}
\maketitle
\begin{abstract}
This paper deals with the asymptotic study of the so-called canard solutions, which arise in the study of real singularly perturbed ODEs. Starting near an attracting branch of the "slow curve", those solutions are crossing a turning point before following for a while a repelling branch of the "slow curve". Assuming that the turning point is degenerate (or non-generic), we apply a correspondence presented in a recent paper. This application needs the definition of a family of functions $\varphi$ that is studied in a first part. Then, we use the correspondence is used to compute the asymptotic expansion in the powers of the small parameter for the canard solution.
\end{abstract}

\section*{Introduction}
We are concerned with the asymptotic study of the so-called \textit{canard solutions} \cite{MR643399}\cite{MR689524}\cite{MR740001} in real singularly perturbed differential equations which allows a degenerate (or non-generic) turning point. More particularly, we apply the formal correspondence presented in \cite{proceeding2} to compute an asymptotic expansion in the powers of the perturbation parameter \cite{MR919406}\cite{MR887738}\cite{MR783844}.\vspace*{3mm}\\
So, we consider equations
$$\varepsilon\dot{u}=\Psi (t,u,\alpha,\varepsilon)$$
where $\Psi$ is $\mathcal{C}^\infty$, $t\in I\subset\R$, $\alpha$ is a real control parameter, $u$ is a real function of the variables $t$ and $\varepsilon$, and $\varepsilon\in]0,\varepsilon_0[$ is a small real parameter which is tending to $0$.\\
By assuming natural and non-natural hypothesis, we restrict our study \cite{papier1} (or see \cite{proceeding1}\cite{these}) to equations
\begin{equation}\label{eq-p}
\eta^{p+1} \dot{u} = (p+1)t^p u+\alpha t^L +S(t,\alpha) + \eta^{p+1} P(t,u,\alpha,\eta)\end{equation}
where $t\in[-t_0,t_0]$, $\alpha\in\R$, $u$ is a real function of the variables $t$ and $\eta$, $\eta=\varepsilon^{1/(p+1)}\in]0,\eta_0[$, $p$ is odd, $L<p$ is even, and where the functions $S$ and $P$ are $\mathcal{C}^\infty$ in their variables. Furthermore, the function $S$ is such that $S(t,0)=0$, and each of its monomial terms have a valuation, with pounds $1$ for $t$ and $p-L+1$ for $\alpha$, strictly greater than $p+1$.\\
The existence of a canard solution in such equations, which is proved in the articles mentionned below, is also proved in \cite{MR1800821}\cite{MR2197445} with geometrical methods.\vspace*{3mm}\\
In what follows, we assume that $p\geq 3$, and thus (\ref{eq-p}) has a degenerate turning point at $t=0$.\vspace*{3mm}\\
For convenience, we adopt the change of variables $\alpha=\eta^{p-L+1}\underline{\alpha}$ and $u=\eta\underline{u}$, which consists of assuming that the parameter (resp. the function) is equal to $\mathcal{O}_{\eta\rightarrow 0}(\eta^{p-L+1})$ (resp. $\mathcal{O}_{\eta\rightarrow 0}(\eta)$).\\
Such assumptions are taken to guarantee \cite{these} the existence and the uniqueness of the canard solution $(\underline{\alpha}_\eta^*,\underline{u}_\eta^*)$.\\
This change of variables transforms (\ref{eq-p}) to the equation
\begin{equation}\label{eq-gen}
\eta^{p+2} \dot{\underline{u}}(t,\eta) = (p+1)\eta t^p \underline{u}(t,\eta)+\eta^{p-L+1}\underline{\alpha} t^L +S(t,\eta^{p-L+1}\underline{\alpha}) + \eta^{p+1} P(t,\eta\underline{u}(t,\eta),\eta^{p-L+1}\underline{\alpha},\eta)\end{equation}
which is the studied equation in what follows.\\
For convenience, the notations $u(t,\eta)$ and $u_\eta(t)$ are assumed to represent the same object.\vspace*{3mm}\\
Moreover, for each fixed couple $(\underline{\beta}_\eta,\underline{v}_\eta)$, the system
$$\left\{\begin{array}{c}\eta^{p+2} \dot{\underline{u}}_\eta(t) = (p+1)\eta t^p \underline{u}_\eta(t)+\eta^{p-L+1}\underline{\alpha}_\eta t^L +S(t,\eta^{p-L+1}\underline{\alpha}_\eta) + \eta^{p+1} P(t,\eta\underline{v}_\eta(t),\eta^{p-L+1}\underline{\beta}_\eta,\eta)\\ \\u(-t_0,\eta)=0=u(t_0,\eta)\end{array}\right.$$
has a solution $(\underline{\alpha}_\eta,\underline{u}_\eta)$ defined all over $[-t_0,t_0]$. So we can define an operator $\Xi_\eta:(\underline{\beta},\underline{v})\mapsto(\underline{\alpha},\underline{u})$, which was proved \cite{papier1} to be a contraction, with Lipschitz constant equals to $\mathcal{O}_{\eta\rightarrow 0}(\eta)$.\\
And thus, the fixed point theorem implies the existence of a canard solution $(\underline{\alpha}_\eta^*,\underline{u}_\eta^*)$ of (\ref{eq-gen}), which can be obtained as the limit of an iterative sequence $((\underline{\alpha}_\eta^{(n)},\underline{u}_\eta^{(n)}))_n$ defined from $\Xi_\eta$.\vspace*{3mm}\\
In a first part, we recall the correspondence, presented in \cite{proceeding2} by giving fundamental assumptions wich leads to a theorem of existence of an asymptotic expansion in the powers of $\eta$ for the canard solution.\\
Then, the second part of this paper is dedicated to the definition and the study of a family of functions, that we call intermediary \cite{MR0435697}.\\
Finally, it is proved that this choice of functions is relevant to apply the correspondence in the linear case. Unfortunately, it encounters some problems with the uniqueness, in the general case.\vspace*{3mm}\\
Nevertheless, we conclude by proving the existence and the uniqueness of an asymptotic expansion in the powers of $\eta$ for the canard solution $(\underline{\alpha}_\eta^*,\underline{u}_\eta^*)$ of the kind
$$\left(\sum_l\underline{a}_l\eta^l,\sum_k\left(\underline{u}_k(t)+\phi_k\left(\frac{t}{\eta}\right)\right)\eta^k\right)$$
where, for all integer $k$, $\underline{a}_k\in\R$, $\widetilde{\underline{u}}_k\in\mathcal{C}^\infty([-t_1,t_1])$ and $\widetilde{\phi}_k$ has an asymptotic expansion with negative powers in its variable $T$ in $+\infty$ and in $-\infty$ that are identical, and tends to $0$ as $T$ tends to $\pm\infty$.\\
Such expansions are uniformly defined in an appreciable neighborhood of the turning point \cite{MR553107}\cite{MR958913}\cite{MR1292646}.\\
Those kind of expansions are currently studied by A. Fruchard and R. Schäfke with methods from the complex analysis.\vspace*{3mm}\\
By convenience, we suppose that the variable $t$ belongs to an interval $[-t_1,t_1]\subset]-t_0,t_0[$. As we are concerned by an asymptotic expansion in an appreciable neighborhood of $0$, this assumption is necessary to avoid boundary problems that are not concerning our study.

\section{Recalls on the correspondence}
In this section, we present the mechanism of the correspondence used to compute the asymptotic expansion in what follows. The fundamental results are demonstrated in \cite{proceeding2} .\vspace*{3mm}\\
To define the spaces which are relevant in this study, we give a finite family $\varphi$ of symbols that have to be chosen in a countable fixed set. All of those symbols are associated to a function, and the family of those functions is denoted by $\widetilde{\varphi}$. The choice of such functions is due to the considered equation.\\
In the remainder of this paper, the symbol $\widetilde{~}$ is put on terms to indicate the substitution of the formal symbols $\varphi$ by their associated topological functions.\vspace*{3mm}\\
Moreover we give an order $\textmd{ord}(.)$, associated to the asymptotic approximation in the powers of $\eta$ to each of the "monomial terms" $t^i\varphi^J\eta^l$.\vspace*{3mm}\\
Those objects allow us to define the sets $\mathcal{A}_k$, as the vector spaces which are generated by terms made up of powers of $t$, $\eta$, and of the intermediary functions $\varphi$, that have an order $\textmd{ord}(.)$ lesser or equals $k$.\\
Then, for all integer $k$, we define
$$\mathcal{B}_k:=\left\{t^i\varphi^J\eta^l:\ |J|\geq 1,\ \textmd{ord}(i,J,l)= k\right\}\textmd{ , and }\mathcal{C}_k:=\R.\eta^k\times\left(\mathcal{C}^\infty([-t_0,t_0]).\eta^k\oplus \textmd{Vect}\mathcal{B}_k\right)$$
$$\mathcal{A}_k:=\mathcal{A}_{k-1}\oplus\mathcal{C}_k=\oplus_{l\leq k}\mathcal{C}_l\textmd{ (where we assume }\mathcal{A}_{-1}:=\mathcal{C}_{-1}:=\{(0,0)\}\textmd{ )}$$
We assume that the order $\textmd{ord}(.)$ is such that, for all $k$, the set $\left\{(i,J,l)\in I;\ |J|\geq 1\textmd{ and }\textmd{ord}(i,J,l)=k\right\}$ is finite. Thus, the substitution in those spaces of the symbols $\varphi$ by the associated functions have sense. In what follows, their images will be respectively denoted by $\widetilde{\mathcal{B}}_k$, $\widetilde{\mathcal{C}}_k$, and $\widetilde{\mathcal{A}}_k$.\\
Finally, for all integer $k$, we denote by $\pi_k$ the natural projection from $\mathcal{A}_{k+1}$ to $\mathcal{A}_k$. The sequence $((\mathcal{A}_k,\pi_k))_k$ has a natural projective limit $(\hat{\mathcal{A}},(\pi_k)_k)$.
\begin{defi}
A sequence $(\Xi_k)_k$, where for all integer $k$ $\Xi_k:\mathcal{A}_k\rightarrow\mathcal{A}_k$, is \textbf{compatible} if, for all $k\in\N$,
$$\pi_k\circ\Xi_{k+1}=\Xi_k\circ\pi_k$$
\end{defi}
\begin{defi}
A sequence of functions $(\Xi_k)_k$ is \textbf{formally equivalent} to a topological operator $\Xi:\ D\subset\mathcal{E}\rightarrow\mathcal{E}$ if and only if, for all $(c_k)_k\in\hat{\mathcal{A}}$, and for all integer $k$:
\begin{itemize}
\item $\widetilde{c}_k$ belongs to the set on which $\Xi$ is well defined.
\item $\Xi(\widetilde{c}_k) -\widetilde{\Xi_k(c_k)}\in\circ_{\eta\rightarrow 0}(\eta^k)$
\end{itemize}
\end{defi}
\begin{defi}
A couple $(\gamma,w)$, where $\gamma\in\R$ and $w:\ (t,\eta)\mapsto w(t,\eta)$ is a $\mathcal{C}^\infty$ function, has the couple $(\sum_k\gamma_k,\sum_kw_k)\in\hat{\mathcal{A}}$ for \textbf{semi-asymptotic expansion} if, for all $n\in\N$, $(\gamma_n,w_n)\in\mathcal{C}_n$, and
$$\forall k\in\N,\ \left|\gamma-\sum_{n=0}^k\widetilde{\gamma}_n\right|=\mathcal{O}_{\eta\rightarrow 0}(\eta^{k+1}),\ \ \sup_{t\in[-t_0,t_0]}\left\{\left|w(t,\eta)-\sum_{n=0}^k\widetilde{w}_n(t,\widetilde{\varphi}(t,\eta),\eta)\right|\right\}=\mathcal{O}_{\eta\rightarrow 0}(\eta^{k+1})$$
The couple $(\gamma,w)$ has an \textbf{asymptotic expansion} if it has an unique semi-asymptotic expansion.
\end{defi}
In practice, the computation of the expansion is a consequence of the two following propositions \cite{proceeding2}:
\begin{prop}\label{exist}
If the choice of intermediary functions $\widetilde{\varphi}$ (and so the construction of the sets $\mathcal{A}_k$) is such that we can define, from the contration $\Xi_\eta$, a sequence of functions $(\Xi_k)_k$, where $\Xi_k$ is defined from $\mathcal{A}_k$ to itself. Moreover, if $(\Xi_k)_k$ is compatible and formally equivalent to $\Xi_\eta$, then we have the \textbf{existence} of a semi-asymptotic expansion for the canard solution $(\underline{\alpha}_\eta^*,\underline{u}_\eta^*)$.
\end{prop}
\noindent More particulary, we consider the coefficients $\underline{u}_k$
$$\underline{u}_k(t,\varphi,\eta)=\check{\underline{u}}_k(t,\varphi)\eta^k$$
\begin{prop}\label{uniq}
If the family $(\mathcal{A}_k)_k$ is \textbf{asymptotically free}, i.e. satisfied to
$$\forall k\in\N\textmd{ fixed, }\forall c_k\in\mathcal{A}_k,\ \left[(\widetilde{c}_k=\mathcal{O}_{\eta\rightarrow 0}(\eta^{k+1}))\ \Rightarrow (c_k=0)\right]$$
Then, we have the \textbf{uniqueness} of the semi-asymptotic expansion that was computed in the previous proposition.
\end{prop}
\noindent For all $k\in\N$, $\Xi_k$ is defined from the substitutions, in the integral form of $\Xi_\eta$, of the functions $P$ and $S$ by their respectives Taylor expansions, truncated at order $k$. And then, by the substitution of all the terms by their associated formal terms in $\mathcal{A}_k$.\\
In practice, for all $k\in\N$, and for all $(\underline{\beta}_k,\underline{v}_k)\in\widetilde{\mathcal{A}}_k$, we denote $(\underline{\alpha}_k,\underline{u}_k):=\Xi_\eta(\underline{\beta}_k,\underline{v}_k)$.\\
As, by definition of $\widetilde{\mathcal{A}}_k$, $(\underline{\beta}_k,\underline{v}_k)$ is associated to $(\hat{\underline{\beta}}_k,\hat{\underline{v}}_k)\in\mathcal{A}_k$, we can define
$$(\hat{\underline{\alpha}}_k,\hat{\underline{u}}_k):=\Xi_k(\hat{\underline{\beta}}_k,\hat{\underline{v}}_k)$$
Finally, the expected asymptotic expansion is obtained as the limit of the sequence $((\hat{\underline{\alpha}}_k,\hat{\underline{u}}_k))_k$.

\section{Study of a family of functions}
In this section, we define the family of intermediary functions which is relevant to our study. This family is defined in two steps:\begin{itemize}
\item[-] Firstly, we define a family of functions $(\psi_k)_k$ from equations that are similar to the studied one.
\item[-] From the properties of those functions, and the definition of their order $\textmd{ord}(.)$, we define the relevant family $(\bar{\psi}_{i,k})_{i,k}$ used for the study, in the linear case, presented in the section 3.
\end{itemize}
Finally, we propose a generalization of this, which brings a theorem in the general case.

\subsection{Definition of a family of functions}
We denote by $\mathcal{I}_\eta(v)$ the canard solution of the system
\begin{equation}\label{eq-interm}
\left\{\begin{array}{c} \eta^{p+1} \dot{w}(t,\eta) = (p+1)t^p w(t,\eta) + \delta_v x^L+ \eta^{p+1} v(t,\eta)\\ \\ \lim_{t\rightarrow -\infty}e^{-(t/\eta)^{p+1}}w(t,\eta)=0=\lim_{t\rightarrow +\infty}e^{-(t/\eta)^{p+1}}w(t,\eta) \end{array} \right.
\end{equation}
where the parameter $\delta_v$ is such that the function $\mathcal{I}_\eta(v)$ is continuous at $t=0$. We have 
$$\delta_v = -\frac{\int_{-\infty}^{+\infty} v(y,\eta)e^{-(y/\eta)^{p+1}}dy}{\int_{-\infty}^{+\infty}y^L e^{-(y/\eta)^{p+1}}dy}\eta^{p+1} = -\frac{\int_{-\infty}^{+\infty} v(y,\eta)e^{-(y/\eta)^{p+1}}dy}{2\int_0^{+\infty}y^L e^{-(y/\eta)^{p+1}}dy}\eta^{p+1}\textmd{ , as $L$ is even}$$
In what follows, we adopt the notation $\lambda_v=\frac{1}{\eta^{p+1}}\delta_v$.\vspace*{3mm}\\
An explicit study of the equation (\ref{eq-interm}) gives:
\begin{prop}\label{def-I}
$$\mathcal{I}_\eta(v)(t)= e^{(t/\eta)^{p+1}}\int_{+\infty}^t\left(v(y,\eta)-\frac{\int_{-\infty}^{+\infty} v(z,\eta)e^{-(z/\eta)^{p+1}}dz}{2\int_0^{+\infty} z^Le^{-(z/\eta)^{p+1}}dz}y^L\right)e^{-(y/\eta)^{p+1}}dy$$
\end{prop}
\noindent In the case $v(t,\eta)=t^k$, we denote
$$\psi_k(T):=\left\{\begin{array}{c}\frac{1}{p+1}\left(\frac{\Gamma\left(\frac{k+1}{p+1}\right)}{\Gamma\left(\frac{L+1}{p+1}\right)}\Gamma\left(\frac{L+1}{p+1};T^{p+1}\right)-\Gamma\left(\frac{k+1}{p+1};T^{p+1}\right)\right)\textmd{ , if $k$ is even}\\ \\ -\frac{1}{p+1}\Gamma\left(\frac{k+1}{p+1};T^{p+1}\right) \textmd{ , if $k$ is odd}\end{array}\right.$$
where $\Gamma(\kappa, .)$ is the incomplete Gamma function $\Gamma(\kappa;T):= e^T\int_T^{+\infty} z^{\kappa-1} e^{-z}dz$. This implies that
$$\mathcal{I}_\eta(X^k)(t)=\eta^{k+1}\psi_k\left(\frac{t}{\eta}\right)$$
where the functions $\psi_k$ are independant from $\eta$. Moreover, those functions are such that:
\begin{prop}$\psi_L=0$ and, for all $k\in\{0,\cdots,p-1\}$,  $\psi_k$ is $\mathcal{C}^\infty$ and limited for all $t\in\R$.
\end{prop}
\noindent Those properties are fundamental to define the order $\textmd{ord}(.)$ of those functions.\vspace*{3mm}\\
As the first definition of the operator $\mathcal{I}_\eta$ implies that $\mathcal{I}_\eta(v)$ is an odd (resp. even) function if $v$ is even (resp. odd), we conclude that $\mathcal{I}_\eta(X^k)$ is well defined for negative values of $t$.\vspace*{3mm}\\
Moreover, the following properties are an immediate consequences of the definition of $\mathcal{I}_\eta$:
\begin{prop}
\begin{itemize}
\item[-]If, for all $\eta$ and $t\in\R,\ v(t,\eta)=\mathcal{O}_{\eta\rightarrow 0}(\eta^K)$, then $\mathcal{I}_\eta(v)(t)=\mathcal{O}_{\eta\rightarrow 0}(\eta^{K+1})$
\item[-] $\mathcal{I}_\eta(u+v)=\mathcal{I}_\eta(u)+\mathcal{I}_\eta(v)$
\item[-] $\mathcal{I}_\eta(g(\eta)\ v)=g(\eta)\ \mathcal{I}_\eta(v)$ 
\item[-] $\mathcal{I}_\eta(X^L)=0$ (et $\delta_{t^L}=-\eta^{p+1}$)
\item[-] $\mathcal{I}_\eta(X^p)=\frac{-1}{p+1}\eta^{p+1}$ (et $\delta_{t^p}=0$)
\end{itemize}
\end{prop}
\noindent Finally, we denote $\mathcal{D}_\eta:\ u \mapsto v$ the operator
$$\mathcal{D}_\eta(u):= \dot{u}-(p+1)\frac{1}{\eta^{p+1}}t^pu$$
The interactions between $\mathcal{D}_\eta$ and $\mathcal{I}_\eta$ brings:
\begin{prop}\label{res-simp}
For all integer $k$, and for all function $u$:
$$\mathcal{I}_\eta(X^k\mathcal{I}_\eta(u))=\frac{1}{k+1}\left(t^{k+1}\mathcal{I}_\eta(u)-\mathcal{I}_\eta(X^{k+1}u)-\lambda_{u}\mathcal{I}_\eta(X^{k+L+1})\right)$$
\begin{equation}\label{eq-simpl}\mathcal{I}_\eta(X^pu)\frac{1}{p+1}\eta^{p+1}(\mathcal{I}_\eta(\dot{u})-u)\end{equation}
In particular
$$\mathcal{I}_\eta(X^{k+p+1})=\frac{1}{p+1}\eta^{p+1}\left((k+1)\mathcal{I}_\eta(X^k)- X^{k+1}\right)$$
\end{prop}
\noindent Moreover, $\mathcal{D}_\eta$ gives:
\begin{prop}\label{psi_k}
The family $\left\{\psi_k:\ k\in\{0,\cdots,p-1\}\verb+\+\{L\}\right\}$ is free.
\end{prop}
\noindent \underline{Proof:}\\
We consider a linear combination
$$\sum_{k=1,\ k\not=L}^{p-1}\iota_k\psi_k\left(\frac{t}{\eta}\right)$$
that we supposed to be equals $0$.\\
As $\mathcal{D}_\eta$ is linear, its definition implies that
$$\mathcal{D}_\eta\left(\psi_k\left(\frac{t}{\eta}\right)\right)=\frac{1}{\eta}\left(\frac{t}{\eta}\right)^k+\frac{\lambda_{t^k}}{\eta^{k-L+1}}\left(\frac{t}{\eta}\right)^L$$
which gives
$$\sum_{k=1,\ k\not=L}^{p-1}\frac{\iota_k}{\eta}\left(\frac{t}{\eta}\right)^k+\left(\sum_{k=1,\ k\not=L}^{p-1}\frac{\lambda_{t^k}}{\eta^{k-L+1}}\iota_k\right)\left(\frac{t}{\eta}\right)^L=0$$
Finally, as the family $T^n=\left(\frac{t}{\eta}\right)^n$ is free, we conclude that the family $\left\{\psi_k:\ k\in\{0,\cdots,p-1\}\verb+\+\{L\}\right\}$ is free too.
$$\fin$$

\subsection{Definition of the relevant family of functions}
To define the relevant family of "intermediary" functions, we study the asymptotic order $\textmd{ord}(.)$ of the monomial terms $t^i\psi_k\eta^l$.\vspace*{3mm}\\
As, for all $m\in]0,1[$ and for all positive real $U$, we have
$$\Gamma(m,U)=e^{-U}\int_U^{+\infty}z^{m-1}e^{-U}dU$$
then, for all $k\in\{0,\cdots,p\}$, we have the following asymptotic approximation of $\psi_k$ near $\pm\infty$:
$$\psi_k(T)\sim_{T\rightarrow \pm\infty}\left\{\begin{array}{c}T^{-(p-\max\{k,L\})}\textmd{ , if $k$ is even}\\ \\ T^{-(p-k)}\textmd{ , if $k$ is odd}\end{array}\right.\sim_{T\rightarrow \pm\infty}T^{-(p-\max\{k,L\})}$$
Moreover, as $t=\mathcal{O}_{\eta\rightarrow 0}(1)$, we have the approximation:
\begin{footnotesize}
$$t^i \psi_k\left(\frac{t}{\eta}\right)\eta^l = \left\{\begin{array}{c}\mathcal{O}_{\eta\rightarrow 0}\left(\eta^{i+l}\psi_k(\mathcal{O}_{\eta\rightarrow 0}(1))\right)=\mathcal{O}_{\eta\rightarrow 0}(\eta^{i+l}) \textmd{ , if } t=\mathcal{O}_{\eta\rightarrow 0}(\eta) \\ \\ \mathcal{O}_{\eta\rightarrow 0}(1)\ \mathcal{O}_{\eta\rightarrow 0}(\left(\frac{t}{\eta}\right)^{-(p-\max\{k,L\})}\eta^l)=\mathcal{O}_{\eta\rightarrow 0}(\eta^{p-\max\{k,L\}+l}) \textmd{ , if } t\not= \mathcal{O}_{\eta\rightarrow 0}(\eta)\end{array}\right.$$
\end{footnotesize}
In conclusion, for all $(i,k,l)\in\N\times\{0,\cdots,p\}\times\N$, the relevant choice of order for the monomial term $t^i \psi_k(\frac{t}{\eta})\eta^l$ is
$$\textmd{ord}(t^ i\psi_k\eta^l):=\min\{p-\max\{k,L\},i\}+l\textmd{ , and }\left(\textmd{ord}(i,k,l):=\right)\textmd{ord}(t^ i\hat{I}_k\eta^l):=\textmd{ord}(t^ i\psi_k\eta^{l+k+1})$$
The choice of order for the term $t^i\eta^l$ is naturally
$$\textmd{ord}(t^i\eta^l):=l$$
This choice implies that the terms of order $K$ are the couples 
\begin{footnotesize}
$$\left(a_K\eta^K,u_K(t)\eta^K+\sum_{(i,k,l);\ i<p-\max\{k,L\},\ i+l=K}u_{i,k}t^i\hat{I}_k\eta^{K-i}+\sum_{(k,l);\ p-\max\{k,L\}+l=K}t^{p-\max\{k,L\}}f_{k,l}(t)\hat{I}_k\eta^l\right)$$
\end{footnotesize}
where $a_K\in\R$, $u_{i,k,l}\in\R$, $u_K$ and $f_{k,l}$ are $\mathcal{C}^\infty$ in $t\in[-t_1,t_1]$.\vspace*{3mm}\\
As there is not a finite number of such couples, we have to define an other family of functions that is using the following terminology:\\
\textit{If $\rho$ is a function of the variable $T$ such that
$$\rho(T)\sim_{T\rightarrow \pm\infty}\sum_{n=0}^N p_nT^n+\sum_{m\geq 1} p_{-m}T^{-m}$$
We denote
$$[\rho(T)]:=\rho(T)-\sum_{n=0}^N p_nT^n$$}
This terminology gives
\begin{prop}
For all $(i,l)\in\N^2$, and $k\in\{0,\cdots,p-1\}\verb+\+\{L\}$, we have
$$\textmd{ord}(t^i\eta^l):=l\textmd{ , and }\textmd{ord}\left(\left[t^i\psi_k\left(\frac{t}{\eta}\right)\right]\eta^l\right)=\textmd{ord}([T^i\psi_k(T)]\eta^{i+l}):=i+l$$
\end{prop}
So, for all integer $K$, our relevant space $\mathcal{C}_K$ is the vector space which contains the couples
$$\left(a_K\eta^K,\ \left(u_K(t)+\sum_{i=0}^K\sum_{k=0,\ k\not=L}^{p-1}u_{i,k}\bar{\psi}_{i,k}\left(\frac{t}{\eta}\right)\right)\eta^K\right)$$
where $a_K\in\R$, $u_{i,k}\in\R$, $u_K\in\mathcal{C}^\infty([-t_1,t_1])$, where we adopt the notation
$$\bar{\psi}_{i,k}(T):=[T^i\psi_k(T)]$$
Thus, for all integer $k$, we define $\mathcal{A}_k:=\bigoplus_{l=0}^k\mathcal{C}_l$.

\section{Computation in the linear case}
This section is dedicated to the computation in the linear case. We will prove that the method and the choice of intermediary functions presented previously allows the computation of the expected asymptotic expansion for the canard solution.\\
In the last section, we will give some comments on the computation of the asymptotic expansions in general. It will be showned that problem of uniqueness of such expansions still remains.\vspace*{3mm}\\
In this section we assume that the function $P$ in the equation (\ref{eq-gen}) is linear in $\underline{u}$.\\
In a first part, using the proposition \ref{exist}, we prove the existence of semi-asymptotic expansions. Then, it is proved, in a second part, that the assumptions of the proposition \ref{uniq} are satisfied.

\subsection{Existence of semi-asymptotic expansions}
To prove the existence of semi-asymptotic expansions for the canard solution, we need:
\begin{lemm}\label{struct-I}For all integer $k$, $\mathcal{I}_\eta$ has a formal equivalent in $\mathcal{A}_k$.
\end{lemm}
\noindent \underline{Proof:}\\
By definition of the spaces $\mathcal{A}_k$, it is sufficient to prove that, for all integer $K$ and for all function $w\in\mathcal{C}^\infty([-t_1,t_1])$, the functions $\mathcal{I}_\eta(w(X))$ and $\mathcal{I}_\eta(\bar{\psi}_{i,k}(X/\eta))\eta^i$ have a natural equivalent formal term in $\mathcal{A}_K$.\\
In what follows, we suppose that an integer $K$ is fixed.\vspace*{3mm}\\
$\otimes$ \textit{Study of terms of the kind $\mathcal{I}_\eta(w(X))$:}\\
We assume that $K>p$ (the case $K\leq p$ is a simplier case).\\
Considering the Taylor expansion of order $p$ for the function $w$, the linearity of $\mathcal{I}_\eta$ gives
$$\mathcal{I}_\eta(w(X))=\sum_{n=0}^p w_{1,n} \mathcal{I}_\eta(X^n)+\mathcal{I}_\eta(X^pR_1(X))$$
Then, the equation(\ref{eq-simpl}) implies that
$$\mathcal{I}_\eta(w(X))=\sum_{n=0}^p w_{1,n} \mathcal{I}_\eta(X^n)-\frac{1}{p+1}\eta^{p+1}R_1+\frac{1}{p+1}\mathcal{I}_\eta(\dot{R}_1(X))\eta^{p+1}$$
By iterating this process, we arrive, for all integer $q$, to an equation of the kind
$$\mathcal{I}_\eta(w(X))=\sum_{n=0}^p \left(\sum_{m=1}^qw_{m,n}\eta^{(m-1)(p+1)}\right) \mathcal{I}_\eta(X^n)-\sum_{m=1}^qS_m\eta^{m(p+1)}+\mathcal{I}_\eta(\Lambda)\eta^{q(p+1)}$$
that we can written
$$\mathcal{I}_\eta(w(X))=\sum_{n=0}^p \left(\sum_{m=1}^qw_{m,n}\eta^{(m-1)(p+1)}\right) \eta^{n+1}[\psi_n(t/\eta)]-\sum_{m=1}^qT_m\eta^{m(p+1)}+\mathcal{I}_\eta(\Lambda)\eta^{q(p+1)}$$
where $w_{m,n}\in\R$, and the functions $T_m$ and $\Lambda$ belongs to $\mathcal{C}^\infty([-t_1,t_1])$.\vspace*{3mm}\\
In conclusion, by denoting $q_K:=\min\{m\in\N;\ m(p+1)>K\}$, the equivalent term of $\mathcal{I}_\eta(w(X))$ which belongs to $\mathcal{A}_K$ is
$$\left\{\begin{array}{c}\sum_{n=0}^{K-1} \left(\sum_{m=1}^qw_{m,n}\eta^{(m-1)(p+1)+n+1}\right)[\psi_n]\textmd{ , if }K\leq p\\ \\ \sum_{n=0}^p \left(\sum_{m=1}^qw_{m,n}\eta^{(m-1)(p+1)+n+1}\right)[\psi_n]-\sum_{m=1}^{q_K-1}T_m(t)\eta^{m(p+1)}\textmd{ , if }K>p\end{array}\right.\in\mathcal{A}_K$$
$\otimes$ \textit{Study of terms of the kind $\mathcal{I}_\eta(\bar{\psi}_{i,k}(X/\eta))\eta^i$:}\\
The formal substitution of a series in negative powers of $T$ in the differential equations on which $\psi_k$ is solution shows that the function $\psi_k$ has an expansion of the kind:
$$\psi_k(T)\sim_{T\rightarrow \pm\infty}\sum_{n\geq p-\max\{k,L\}}\rho_{-n}T^{-n}$$
So, for all $i\in\N$, and for all $k\in\{0,\cdots,p-1\}\verb+\+\{L\}$, we have:
\begin{equation}\label{croc-psi}
\bar{\psi}_{i,k}\left(\frac{t}{\eta}\right)=\left\{\begin{array}{c} \left(\frac{t}{\eta}\right)^i\psi_k\left(\frac{t}{\eta}\right)\textmd{ , if } i<p-\max\{k,L\}\\ \\ \left(\frac{t}{\eta}\right)^i\psi_k\left(\frac{t}{\eta}\right)-\sum_{n=0}^{i-p+\max\{k,L\}}\rho_{-n}\left(\frac{t}{\eta}\right)^{i-n}\textmd{ , if } i\geq p-\max\{k,L\}\end{array}\right.
\end{equation}
In what follows, we study the second case, as the first one can be prove with a similar method.\\
So, we assume that $i\geq p-\max\{k,L\}$, and then
$$\mathcal{I}_\eta\left(\bar{\psi}_{i,k}\left(\frac{X}{\eta}\right)\right)\eta^i=\mathcal{I}_\eta\left(X^i\psi_k\left(\frac{X}{\eta}\right)\right)-\sum_{n=0}^{i-p+\max\{k,L\}}\rho_{-n}\mathcal{I}_\eta(X^{i-n})\eta^n$$
The proposition \ref{res-simp} brings
\begin{footnotesize}
$$\mathcal{I}_\eta\left(X^i\psi_k\left(\frac{X}{\eta}\right)\right)=\frac{1}{\eta^{k+1}}\mathcal{I}_\eta\left(X^i\mathcal{I}_\eta(X^k)\right)=\frac{1}{(i+1)\eta^{k+1}}\left(t^{i+1}\mathcal{I}_\eta(X^k)-\mathcal{I}_\eta(X^{i+k+1})-\lambda_{X^k}\mathcal{I}_\eta(X^{i+L+1})\right)$$
\end{footnotesize}
which gives
\begin{small}
$$\mathcal{I}_\eta\left(X^i\psi_k\left(\frac{X}{\eta}\right)\right)=\frac{1}{i+1}\left(t^{i+1}\psi_k\left(\frac{t}{\eta}\right)-\eta^{i+1}\psi_{i+k+1}\left(\frac{t}{\eta}\right)-\frac{\lambda_{t^k}}{\eta^{k-L}}\eta^{i+1}\psi_{i+L+1}\left(\frac{t}{\eta}\right)\right)=\ldots$$
$$=\frac{1}{i+1}\eta^{i+1}\left(\left(\frac{t}{\eta}\right)^{i+1}\psi_k\left(\frac{t}{\eta}\right)-\psi_{i+k+1}\left(\frac{t}{\eta}\right)-\frac{\lambda_{t^k}}{\eta^{k-L}}\psi_{i+L+1}\left(\frac{t}{\eta}\right)\right)$$
\end{small}
In conclusion, the formula (\ref{croc-psi}) implies that $\mathcal{I}_\eta\left(X^i\psi_k\left(\frac{X}{\eta}\right)\right)$ belongs to $\bigcup_k\mathcal{A}_k$.\\
Even if it means to make truncations, we conclude that $\mathcal{I}_\eta(\bar{\psi}_{i,k}(X/\eta))\eta^i$ has an equivalent formal term in $\mathcal{A}_K$.
$$\fin$$
Thus, we can prove the fundamental result of this section:
\begin{thme}\label{exist-lin}
There exists a sequence $(\Xi_K)_K$, which is defined from the topological contraction $\Xi_\eta$, which is formally equivalent to $\Xi_\eta$ and compatible.
\end{thme}
\noindent In particular, this theorem allows the computation of semi-asymptotic expansions for the canard solution.\vspace*{3mm}\\
\underline{Demonstration of the theorem:}\\
This sequence is defined by iteration. So, an integer $K$ is supposed to be fixed.\\
To define $\Xi_K$, a couple $(\underline{\beta}_K,\underline{v}_K)\in\widetilde{\mathcal{A}}_K$ is given. We recall that its equivalent formal term is denoted by $(\hat{\underline{\beta}}_K,\hat{\underline{v}}_K)\in\mathcal{A}_K$. We denote $(\underline{\alpha}_K,\underline{u}_K):=\Xi_\eta(\underline{\beta}_K,\underline{v}_K)$.\vspace*{3mm}\\
\textit{Study of the equivalent formal term of $\underline{\alpha}_K$:}\vspace*{3mm}\\
The definition of the operator $\Xi_\eta$ implies that the parameter $\underline{\alpha}_K$ satisfies to
$$0=\eta^{p-L}\underline{\alpha}_K\ \int_{-t_0}^{t_0}\xi^L e^{-(\xi/\eta)^{p+1}} d\xi+\ldots$$
$$+\frac{1}{\eta}\int_{-t_0}^{t_0} S(\xi,\eta^{p-L+1}\underline{\alpha}_K)e^{-(\xi/\eta)^{p+1}} d\xi+\eta^p\int_{-t_0}^{t_0}P(\xi,\eta\underline{v}_K(\xi,\eta),\eta^{p-L+1}\underline{\beta}_K,\eta)e^{-(\xi/\eta)^{p+1}} d\xi$$
By substituting $S$ by its explicit formulation and by substituting to the $\mathcal{C}^\infty$-function $P$ the term of order $K$ in its Taylor expansion, the right-handed term of this equation become
\begin{footnotesize}
$$\eta^{p-L}\underline{\alpha}_K\ \int_{-t_0}^{t_0}\xi^L e^{-(\xi/\eta)^{p+1}} d\xi+\ldots$$
$$+\frac{1}{\eta}\left(\sum_{i=0}^{p+1}\eta^{(j_i+1)(p-L+1)}\underline{\alpha}_K^{j_i+1}S_i(\eta^{p-L+1}\underline{\alpha}_K)\int_{-t_0}^{t_0} \xi^i e^{-(\xi/\eta)^{p+1}} d\xi+\right.\ldots$$
$$\left.+\int_{-t_0}^{t_0}\xi^{p+2}R(\xi,\eta^{p-L+1}\underline{\alpha}_K)e^{-(\xi/\eta)^{p+1}} d\xi\right)+\ldots$$
$$+\eta^p\sum_{c\in\{0,1\}}\sum_{a,b,d=0}^{K+1}p_{a,b,c,d}\eta^{b+c(p-L+1)+d}\underline{\beta}_K^c\int_{-t_0}^{t_0}\xi^a\underline{v}_K(\xi,\eta)^b\ e^{-(\xi/\eta)^{p+1}} d\xi$$
\end{footnotesize}
which is supposed, by approximation, to be $0$.\\
As $(\underline{\beta}_K,\underline{v}_K)$ belongs to $\widetilde{\mathcal{A}}_K$, we can replace the terms $\underline{\beta}_K$ et $\underline{v}_K$ by the topological version of $(\hat{\underline{\beta}}_K,\hat{\underline{v}}_K)\in\mathcal{A}_K$.\\
Finally, an approximation of the integral terms brings to an equation of the kind
$$E(\underline{\alpha}_K,\eta)=0 \textmd{ , where $E$ is } \mathcal{C}^\infty\textmd{ in its variables}$$
More precisely, an asymptotic approximation of each integral terms show that all of the monomial terms which constituate the function $E$, and contains $\underline{\alpha}_K$, are asymptotically dominated by the term $\eta^{p-L}\underline{\alpha}_K\ \int_{-t_0}^{t_0}\xi^L e^{-(\xi/\eta)^{p+1}} d\xi$.\\
Applying the implicit function theorem, we conclude that the parameter $\underline{\alpha}_K$ is a $\mathcal{C}^\infty$-function of the variable $\eta$. And, we denote by $\hat{\underline{\alpha}}_K$ the term of order $K$ in its Taylor expansion.\vspace*{3mm}\\
\textit{Study of the equivalent formal term of  $\underline{u}_K$:}\vspace*{3mm}\\
By definition of the operator $\Xi_\eta$, the function $\underline{u}_K$ satisfies, for all $t$ and $\eta$, to
$$\underline{u}_K(t,\eta)=\frac{1}{\eta^{p+1}}e^{(t/\eta)^{p+1}}\int_{t_0}^t\left(\eta^{p-L}\underline{\alpha}_K\xi^L+\right.\ldots$$
$$\left.+\frac{1}{\eta} S(\xi,\eta^{p-L+1}\underline{\alpha}_K) + \eta^p P(\xi,\eta\underline{v}_K(\xi,\eta),\eta^{p-L+1}\underline{\beta}_K,\eta)\right) e^{-(\xi/\eta)^{p+1}} d\xi$$
The definition of the functions $S$ and $P$ show that the presence of the terms $\frac{1}{\eta^{p+1}}$ and $\frac{1}{\eta}$ did not prevent each of the terms to be equals $\mathcal{O}_{\eta\rightarrow 0}(1)$.\vspace*{3mm}\\
In this last expression, we are willing to replace $\underline{\alpha}_K$ by a formula which comes from
$$0=\eta^{p-L}\underline{\alpha}_K\ \int_{-t_0}^{t_0}\xi^L e^{-(\xi/\eta)^{p+1}} d\xi+\ldots$$
$$+\frac{1}{\eta}\int_{-t_0}^{t_0} S(\xi,\eta^{p-L+1}\underline{\alpha}_K)e^{-(\xi/\eta)^{p+1}} d\xi+\eta^p\int_{-t_0}^{t_0}P(\xi,\eta\underline{v}_K(\xi,\eta),\eta^{p-L+1}\underline{\beta}_K,\eta)e^{-(\xi/\eta)^{p+1}} d\xi$$
So, we remark that the term
\begin{small}\begin{equation}\label{eq-a}
\eta^{p-L}\underline{\alpha}_K\ \int_{-\infty}^{+\infty}\xi^L e^{-(\xi/\eta)^{p+1}} d\xi+\ldots$$
$$+\frac{1}{\eta}\int_{-\infty}^{+\infty} S(\xi,\eta^{p-L+1}\underline{\alpha}_K)e^{-(\xi/\eta)^{p+1}} d\xi+\eta^p\int_{-\infty}^{+\infty}P(\xi,\eta\underline{v}_K(\xi,\eta),\eta^{p-L+1}\underline{\beta}_K,\eta)e^{-(\xi/\eta)^{p+1}} d\xi\end{equation}\end{small}
is exactly
\begin{footnotesize}
$$\eta^{p-L}\underline{\alpha}_K\ \int_{-\infty}^{-t_0}\xi^L e^{-(\xi/\eta)^{p+1}} d\xi +\ldots$$
$$+\frac{1}{\eta}\int_{-\infty}^{-t_0} S(\xi,\eta^{p-L+1}\underline{\alpha}_K)e^{-(\xi/\eta)^{p+1}} d\xi +\eta^p\int_{-\infty}^{-t_0}P(\xi,\eta\underline{v}_K(\xi,\eta),\eta^{p-L+1}\underline{\beta}_K,\eta)e^{-(\xi/\eta)^{p+1}} d\xi-\ldots$$
$$-\eta^{p-L}\underline{\alpha}_K\ \int_{t_0}^{+\infty}\xi^L e^{-(\xi/\eta)^{p+1}} d\xi -\ldots$$
$$-\frac{1}{\eta}\int_{t_0}^{+\infty} S(\xi,\eta^{p-L+1}\underline{\alpha}_K)e^{-(\xi/\eta)^{p+1}} d\xi -\eta^p\int_{t_0}^{+\infty}P(\xi,\eta\underline{v}_K(\xi,\eta),\eta^{p-L+1}\underline{\beta}_K,\eta)e^{-(\xi/\eta)^{p+1}} d\xi$$
\end{footnotesize}
By the change of variable $\xi \pm t_0=\zeta$, each terms can be written $e^{-(t_0/\eta)^{p+1}}.f(\underline{\alpha}_K,\eta)$, where the function $f$ is $\mathcal{C}^\infty$. Consequently, each of those six terms is exponentially small, and so can be put on side of our study, as a term which has a regular asymptotic expansion in the powers of $\eta$ equals $0$. In conclusion, we approximate in what follows the term (\ref{eq-a}) by $0$.\vspace*{3mm}\\
This approximation allows us to substitute the term $\eta^{p-L}\underline{\alpha}_K$, in the formulation of $\underline{u}_K$ previously given, by 
$$-\frac{1}{\eta}\frac{\int_{-\infty}^{+\infty} S(s,\eta^{p-L+1}\underline{\alpha}_K)e^{-(s/\eta)^{p+1}} ds}{\int_{-\infty}^{+\infty}s^L e^{-(s/\eta)^{p+1}} ds}-\eta^p\frac{\int_{-\infty}^{+\infty}P(s,\eta\underline{v}_K(s,\eta),\eta^{p-L+1}\underline{\beta}_K,\eta)e^{-(s/\eta)^{p+1}} ds}{\int_{-\infty}^{+\infty}s^L e^{-(s/\eta)^{p+1}} ds}$$
This substitution gives for approximation of $\underline{u}_K$:
\begin{footnotesize}
$$\frac{1}{\eta^{p+1}}e^{(t/\eta)^{p+1}}\left[\int_{t_0}^t\frac{1}{\eta} \left(S(\xi,\eta^{p-L+1}\underline{\alpha}_K)-\frac{\int_{-\infty}^{+\infty} S(s,\eta^{p-L+1}\underline{\alpha}_K)e^{-(s/\eta)^{p+1}} ds}{2\int_0^{+\infty}s^L e^{-(s/\eta)^{p+1}} ds}\xi^L \right) e^{-(\xi/\eta)^{p+1}} d\xi+\right.\ldots$$
$$+\left.\eta^p\int_{t_0}^t \left(P(\xi,\eta\underline{v}_K(\xi,\eta),\eta^{p-L+1}\underline{\beta}_K,\eta)-\frac{\int_{-\infty}^{+\infty} P(s,\eta\underline{v}_K(s,\eta),\eta^{p-L+1}\underline{\beta}_K,\eta)e^{-(\frac{s}{\eta})^{p+1}} ds}{2\int_0^{+\infty}s^L e^{-(s/\eta)^{p+1}} ds}\xi^L \right) e^{-(\frac{\xi}{\eta})^{p+1}} d\xi\right]$$
\end{footnotesize}
As $t\in[-t_1,t_1]\subset]-t_0,y_0[$, we make the approximation which consists in replacing $t_0$ by $+\infty$ in the integral terms. This substitution generate only terms which have an asymptotic approximation equal to $0$.\\
Moreover, an asymptotic approximation of each integrals show that all of the terms which defined $\underline{u}_K$ are equal $\mathcal{O}_{\eta\rightarrow 0}(1)$. And so, in what follows, we put away the negative powers of $\eta$ that appear.\vspace*{3mm}\\
To propose an equivalent formal term of $\underline{u}_K$, the first step is to substitute formally the functions $S$ and $P$ by the term of order $K$ of their respective Taylor expansions. Then, to substitute to $\underline{\beta}_K$ (resp. $\underline{\alpha}_K$, $\underline{v}_K$) the topological version of $\hat{\underline{\beta}}_K$ (resp. $\hat{\underline{\alpha}}_K$, $\hat{\underline{v}}_K$). This gives
\begin{footnotesize}
$$e^{(t/\eta)^{p+1}}\int_{+\infty}^t\left( w_n(\xi)\eta^n-\frac{\int_{-\infty}^{+\infty} w_n(s)\eta^ne^{-(s/\eta)^{p+1}} ds}{\int_{-\infty}^{+\infty}s^L e^{-(s/\eta)^{p+1}} ds}\xi^L \right) e^{-(\xi/\eta)^{p+1}} d\xi\textmd{ , où }w_n\in\mathcal{C}^\infty([-t_1,t_1])$$
$$\textmd{or }e^{(t/\eta)^{p+1}}\int_{+\infty}^t\left( \bar{\psi}_{n,k}(\xi/\eta)\eta^m-\frac{\int_{-\infty}^{+\infty} \bar{\psi}_{n,k}(s/\eta)\eta^me^{-(s/\eta)^{p+1}} ds}{\int_{-\infty}^{+\infty}s^L e^{-(s/\eta)^{p+1}} ds}\xi^L \right) e^{-(\xi/\eta)^{p+1}} d\xi$$
\end{footnotesize}
where $n$ and $m$ are integer, and $k\in\{0,\cdots,p-1\}\verb+\+\{L\}$.\vspace*{3mm}\\
The proposition \ref{def-I} shows that those terms are respectively $\mathcal{I}_\eta(w_n(X)\eta^n)$ and $\mathcal{I}_\eta(X^n\mathcal{I}_\eta(X^k)(X)\eta^m)$.\\
Moreover, even if it means to took away terms that have an order $\textmd{ord}(.)$ greater than $K+1$, the lemma \ref{struct-I} implies that we are able to define a term $\hat{\underline{u}}_K$, which is formally equivalent to $\underline{u}_K$ in $\mathcal{A}_K$, such that
$$(\hat{\underline{\alpha}}_K,\hat{\underline{u}}_K)\in\mathcal{A}_K$$
Consequently, the operator $\Xi_K$ can be defined by $\Xi_K(\hat{\underline{\beta}}_K,\hat{\underline{v}}_K):=(\hat{\underline{\alpha}}_K,\hat{\underline{u}}_K)$.\vspace*{3mm}\\
In conclusion, this construction shows that $(\Xi_k)_k$ is formally equivalent to $\Xi_\eta$, and is compatible.
$$\fin$$
In the next part, we study the uniqueness of the semi-asymptotic expansion for the canard solution that we can compute.

\subsection{Uniqueness of the asymptotic expansion}
This part is dedicated to the proof of the following theorem:
\begin{thme}The sequence $(\mathcal{A}_k)_k$ is asymptotically free, i.e.
$$\forall K\in\N\textmd{ fixed, }\forall c_K\in\mathcal{A}_K,\ \left[\ \left(\widetilde{c}_K=\mathcal{O}_{\eta\rightarrow 0}(\eta^{K+1})\right)\Rightarrow c_K=0\ \right]$$
\end{thme}
\noindent Its demonstration needs:
\begin{lemm}
For each fixed integer $i$, the family
$$\mathcal{E}_i:=\left\{T^n:\ n\in\{0,\cdots,p\}\}\cup\{T^i\psi_k(T):\ k\in\{0,\cdots,p-1\}\verb+\+\{L\}\right\}$$
is free
\end{lemm}
\noindent \underline{Proof of the lemma:}\\
This proof use the definition and the properties of the operator $\mathcal{D}_\eta$ described above. In particular, we note that
\begin{footnotesize}
$$\mathcal{D}_\eta(1)=-\frac{p+1}{\eta}\left(\frac{t}{\eta}\right)^p\textmd{ , }\forall n\in\{1,\cdots,p\},\ \mathcal{D}_\eta\left(\left(\frac{X}{\eta}\right)^n\right)=n\eta\left(\frac{t}{\eta}\right)^{n-1}-(p+1)\eta^{2n-1}\left(\frac{t}{\eta}\right)^{p+n}\ ,$$
$$\forall k\in\{0,\cdots,p-1\}\verb+\+\{L\},\ \mathcal{D}_\eta\left(\psi_k\left(\frac{t}{\eta}\right)\right)=\frac{1}{\eta}\left(\frac{t}{\eta}\right)^k+\frac{\lambda_{t^k}}{\eta^{k-L+1}}\left(\frac{t}{\eta}\right)^L\textmd{ , et }$$
$$\forall i\in\N^*,\ \mathcal{D}_\eta\left(\left(\frac{t}{\eta}\right)^i\psi_k\left(\frac{t}{\eta}\right)\right)=\frac{1}{\eta}\left(\frac{t}{\eta}\right)^{i+k}+\frac{\lambda_{t^k}}{\eta^{k-L+1}}\left(\frac{t}{\eta}\right)^{i+L}+\frac{i}{\eta}\left(\frac{t}{\eta}\right)^{i-1}\psi_k\left(\frac{t}{\eta}\right)$$
\end{footnotesize}
The demonstration consists in a recurrence over the integer $i$:\vspace*{3mm}\\
\underline{\textit{Initialization:}} By definition, $$\mathcal{E}_0:=\left\{T^n:\ n\in\{0,\cdots,p\}\}\cup\{\psi_k(T):\ k\in\{0,\cdots,p-1\}\verb+\+\{L\}\right\}$$
So, we supposed that 
$$\mu_0+\sum_{n=1}^p\mu_nT^n+\sum_{k=0,\ k\not=L}^{p-1}\nu_k\psi_k(T)=0$$
The application of the linear operator $\mathcal{D}_\eta$ transforms this equation to
\begin{footnotesize}
$$-\frac{p+1}{\eta}\mu_0T^p+\sum_{n=1}^p\left(n\eta\mu_n T^{n-1}-(p+1)\eta^{2n-1}\mu_nT^{p+n}\right)
+\sum_{k=0,\ k\not=L}^{p-1}\left(\frac{1}{\eta}\nu_kT^k+\frac{\lambda_{t^k}}{\eta^{k-L+1}}\nu_kT^L\right)=0$$
\end{footnotesize}
Moreover, as the family $\{T^m:\ m\in\{0,\cdots,2p\}\}$ is free and as, for all $n\in\{0,\cdots,p\}$, the term $-(p+1)\eta^{2n-1}\mu_nT^{p+n}$ is the only one which allows $T^{p+n}$ in factor, we conclude that all the coefficients $\mu_n$ are zero.\\
Finally, as $k\not=L$, we have $\nu_k=0$.\vspace*{3mm}\\
\underline{\textit{Heredity:}} We fixed an integer $i\geq 1$.\\
The family $\mathcal{D}_\eta(\mathcal{E}_i)$ is constituated by the terms:
$$-\frac{p+1}{\eta}T^p,\ n\eta T^{n-1}-(p+1)\eta^{2n-1}T^{p+n},\ \frac{1}{\eta}T^{i+k}+\frac{\lambda_{t^k}}{\eta^{k-L+1}}T^{i+L}+\frac{i}{\eta}T^{i-1}\psi_k(T)$$
where $n\in\{0,\cdots,p\}$, and $k\in\{0,\cdots,p-1\}\verb+\+\{L\}$.\vspace*{3mm}\\
All of those terms are generated by $\mathcal{E}_{i-1}$ which is, by hypothesis of recurrence, free.\\
We consider, at first, the terms which contains $T^{p+n}$. By application of the proposition \ref{psi_k}, the family of the functions $\psi_k$ is free. So, we conclude that the family $\mathcal{E}_i$ is free.
$$\fin$$
\underline{Demonstration of the theorem:}\\
By definition of $\mathcal{A}_K$, each $c_K\in\mathcal{A}_K$ is such that
$$c_K=\left(\sum_{n=0}^Ka_n\eta^n,\sum_{n=0}^K\left(u_n(t)+\sum_{i=0}^n\sum_{k=0,\ k\not=L}^{p-1}u_{i,k}\bar{\psi}_{i,k}\right)\eta^n\right)$$
The uniqueness of the coefficients of the first member being immediate, we will be concerned by the second member. The proof consists in a recurrence over the integer $K$.\vspace*{3mm}\\
An element $c_K\in\mathcal{A}_K$ is supposed to be given, and supposed to be such that $\widetilde{c}_K=\mathcal{O}_{\eta\rightarrow 0}(\eta^{K+1})$. By definition we have, for all $t\in[-t_1,t_1]$:
\begin{equation}\label{B}
\left|\sum_{n=0}^K\left(u_n(t)+\sum_{i=0}^n\sum_{k=0,\ k\not=L}^{p-1}u_{i,k}\left[\left(\frac{t}{\eta}\right)^i\widetilde{\psi}_k\left(\frac{t}{\eta}\right)\right]\right)\eta^n\right|=\mathcal{O}_{\eta\rightarrow 0}(\eta^{K+1})
\end{equation}
To prove the initialization ($K=0$), the definition of $\textmd{ord}(.)$ gives
$$\forall t\in[-t_1,t_1],\ \left|u_0(t)+\sum_{k=0,\ k\not=L}^{p-1}u_{0,k}\left[\widetilde{\psi}_k\left(\frac{t}{\eta}\right)\right]\right|=\mathcal{O}_{\eta\rightarrow 0}(\eta)$$
A similar demonstration than the one proposed below to prove the heredity shows that the function $u_0$ is $0$, and that $u_{0,k}=0$.\vspace*{3mm}\\
To prove the heredity of the recurrence we assume, by hypothesis of recurrence, that the property is true for all integer smaller or equal to $K-1$.\\
As, by definition of $\textmd{ord}(.)$ and $\mathcal{C}_K$, we know that
$$\left|\left(u_K(t)+\sum_{i=0}^K\sum_{k=0,\ k\not=L}^{p-1}u_{i,k}\left[\left(\frac{t}{\eta}\right)^i\widetilde{\psi}_k\left(\frac{t}{\eta}\right)\right]\right)\eta^K\right|=\mathcal{O}_{\eta\rightarrow 0}(\eta^{K})$$
we conclude from (\ref{B}) that
$$\left|\sum_{n=0}^{K-1}\left(u_n(t)+\sum_{i=0}^n\sum_{k=0,\ k\not=L}^{p-1}u_{i,k}\left[\left(\frac{t}{\eta}\right)^i\widetilde{\psi}_k\left(\frac{t}{\eta}\right)\right]\right)\eta^n\right|=\mathcal{O}_{\eta\rightarrow 0}(\eta^K)$$
And so, the hypothesis of recurrence gives
$$\forall n\in\{0,\cdots,K-1\},\ \forall i\in\{0,\cdots,n\},\ \forall k\in\{0,\cdots,p-1\}\verb+\+\{L\},\ u_n=0\textmd{ , and }u_{i,k}=0$$
Finally, (\ref{B}) become
$$\forall t\in[-t_1,t_1],\ \left|u_K(t)+\sum_{i=0}^K\sum_{k=0,\ k\not=L}^{p-1}u_{i,k}\left[\left(\frac{t}{\eta}\right)^i\widetilde{\psi}_k\left(\frac{t}{\eta}\right)\right]\right|=\mathcal{O}_{\eta\rightarrow 0}(\eta)$$
As, by construction, the functions $T\mapsto\left[T^i\widetilde{\psi}_k(T)\right]$ are $\mathcal{C}^\infty$ and limited over $\R$, a study with $t=\mathcal{O}_{\eta\rightarrow 0}(\eta)$ gives $u_K=0$.\\
Moreover, the formula (\ref{croc-psi}) and the previous lemma shows that
$$\left\{\left[T^i\widetilde{\psi}_k(T)\right]:\ i\in\{0,\cdots,K\},\ k\in\{0,\cdots,p-1\}\verb+\+\{L\}\right\}$$
is free.
$$\fin$$
The two theorems presented in this section gives a correspondence that brings to the computation of an asymptotic expansion in the powers of $\eta$ for the canard solution of (\ref{eq-gen}), in the linear case. This can be presented in the more general result:
\begin{thme}
If the equation (\ref{eq-gen}) is linear, its canard solution $(\underline{\alpha}^*,\underline{u}^*)$ has an asymptotic expansion that can be written
$$\left(\sum_na_n\eta^n,\sum_n\left(u_n(t)+\sum_{i=0}^n\sum_{k=0,\ k\not=L}^{p-1}u_{i,k}\bar{\psi}_{i,k}\left(\frac{t}{\eta}\right)\right)\eta^n\right)$$
\end{thme}
\noindent This method gives a formal process which compute this asymptotic expansion in the powers of $\eta$, which remains uniformly valid all over $[-t_1,t_1]$.\vspace*{3mm}\\
\underline{Remark:}\\
Writing the expansion of the function $\underline{u}^*$ on the form
$$\sum_l\left(u_l(t)+\sum_{k=0,\ k\not=L}^{p-1}\phi_{k,l}\left(\frac{t}{\eta}\right)\right)\eta^l$$
gives an asymptotic expansion of the regularly perturbed differential equations which is associated to (\ref{eq-gen}) (obtained by the change of variables $t=\eta T$, $u=\eta U$, $\alpha=\eta^{p-L+1}A$). This last equation is defined in a domain which is growing as $\eta$ tends to $0$.\\
And so, the change of variable $\tau=\frac{t}{\eta}$, brings to a local study of the canard solution in a small neighborhood of $t=0$.

\section{Study in the general case}
In the general case, the family of intermediary functions which was used in the linear case is clearly not sufficient. As the function $P$ is not necessary linear in $\underline{u}$, we have to take into account the functions generated by the product of others intermediary functions, as $\mathcal{I}_\eta(X^k)\mathcal{I}_\eta(X^l)$, and the one that are consequently obtained by composition, as $\mathcal{I}_\eta(X^i\mathcal{I}_\eta(X^k)\mathcal{I}_\eta(X^l))$ for example. Due to the lack of control of those functions by a theorem similar to the proposition \ref{res-simp} for the linear case, those problems still persists.\vspace*{3mm}\\
The natural space of the intermediary functions $\mathcal{F}$ can be defined as the inductive limit of the sequence $(\mathcal{F}_k)_k$ defined as follow:
$$\mathcal{F}_0:=\{1\}\textmd{ et }\forall k\geq 1,\ \mathcal{F}_k:=(t\ \mathcal{F}_{k-1})\cup (\eta\ \mathcal{F}_{k-1})\cup_{k=0}^{p-1} (\mathcal{I}_\eta(X^k)\ \mathcal{F}_{k-1})\cup \mathcal{I}_\eta(\mathcal{F}_{k-1})$$
On those spaces, the natural order is given by
$$\textmd{ord}(t^i\eta^l)=l\textmd{ , and }\textmd{ord}([T^i\psi_{k_1}(T)\ldots\psi_{k_n}(T)]\eta^{i+l})=i+l$$
And it is possible to define the spaces $\mathcal{A}_k$ as the vector spaces generated by the monomial terms that have an order lesser than $k$.\\
Unfortunately, those spaces $\mathcal{A}_k$ are infinite, and although we have determined some properties that allow us to extract a family of functions from $\mathcal{F}$, it has not be proved that this family is such that the expansions that are computed are unique.\vspace*{3mm}\\
On the other hand, even with the infinite family $\mathcal{F}$, a demonstration similar to the theorem \ref{exist-lin} brings:
\begin{prop}The canard solution of (\ref{eq-gen}) has at least one semi-asymptotic expansion in the powers of $\eta$, which is belonging to the space $\hat{\mathcal{A}}$ obtained as a projective limit of $(\mathcal{A}_k)_k$.
\end{prop}
Nevertheless, it is possible to use this result to prove the existence and the uniqueness of a more general form of an asymptotic expansion for the canard solution.\\
This definition needs to consider the set $\mathcal{D}$ which contains the functions of the variable $T$ that have an asymptotic expansion in negative powers of $T$ at $+\infty$ and at $-\infty$ which are identical.\vspace*{3mm}\\
In particular, for all $i\in\N$ and $k\in\{0,\cdots,p-1\}\verb+\+\{L\}$, $\bar{\psi}_{i,k}\in\mathcal{D}$.\vspace*{3mm}\\
Considering this set, the previous proposition gives the existence of a semi-asymptotic expansion for the canard solution of the kind
$$\underline{u}^*(t,\eta)\sim \sum_k \left(\underline{u}_k(t)+\phi_k\left(\frac{t}{\eta}\right)\right)\eta^k$$
where, for all integer $k$, $\underline{u}_k\in\mathcal{C}^\infty([-t_1,t_1])$ and $\phi_k\in \mathcal{D}$ tends to $0$ as $T$ tends to $\pm\infty$.\\
The following theorem gives the uniqueness of this kind of asymptotic expansion:
\begin{thme}
The canard solution $(\underline{\alpha}^*,\underline{u}^*)$ of (\ref{eq-gen}) has an asymptotic expansion in the powers of $\eta$ of the kind
$$\left(\sum_l\underline{a}_l\eta^l,\sum_k\left(\underline{u}_k(t)+\phi_k\left(\frac{t}{\eta}\right)\right)\eta^k\right)$$
where, for all integer $k$, $\widetilde{\underline{u}}_k\in\mathcal{C}^\infty([-t_1,t_1])$ and $\widetilde{\phi}_k\in\mathcal{D}$ and tends to $0$ in $\pm\infty$.
\end{thme}
\noindent \underline{Demonstration:}\\
As the proof of the uniqueness of the expansion for the parameter is trivial, we will concentrate ourselves on the function.\\
Similarly than the proof of the uniqueness of the asymptotic expansion proposed in the linear case we suppose, by hypothesis of recurrence, and for a fixed $K$, that if we have
$$\forall t\in[-t_1,t_1],\ \sum_{k=0}^{K-1}\left(\underline{u}_k(t)+\phi_k\left(\frac{t}{\eta}\right)\right)\eta^k=\mathcal{O}_{\eta\rightarrow 0}(\eta^K)$$
then all of the coefficients of this sum are $0$.\\
We suppose that
$$\forall t\in[-t_1,t_1],\ \sum_{k=0}^K\left(\underline{u}_k(t)+\phi_k\left(\frac{t}{\eta}\right)\right)\eta^k=\mathcal{O}_{\eta\rightarrow 0}(\eta^{K+1})$$
As the term $\left(\underline{u}_K(t)+\phi_K\left(\frac{t}{\eta}\right)\right)\eta^K$ is, by definition, equal to $\mathcal{O}_{\eta\rightarrow 0}(\eta^K)$, the hypothesis of recurrence gives
$$\left(\underline{u}_K(t)+\phi_K\left(\frac{t}{\eta}\right)\right)\eta^K=\mathcal{O}_{\eta\rightarrow 0}(\eta^{K+1})$$
which we write
$$\underline{u}_K(t)+\phi_K\left(\frac{t}{\eta}\right)=\mathcal{O}_{\eta\rightarrow 0}(\eta)$$
Assuming than $\eta$ tends to $0$, the hypothesis on $\phi_K$ implies that $\underline{u}_K=0$.\\
And then, by considering the expansion of the function $\phi_K$, we deduce that $\phi_K=0$.
$$\fin$$
In this result, we have lost the possibility of having an algorithm that computes the expected asymptotic expansion. We have prove that the computed expansion is a semi-asymptotic expansion for the canard solution.\vspace*{3mm}\\
Nevertheless, this theorem remains important because it shows the existence and the uniqueness of an asymptotic expansion in the powers of $\eta$, which remain uniformly valid all over the interval $[-t_1,t_1]$.\vspace*{3mm}\\
\underline{Remark 1:}\\
The change of variables $t=\eta T$, $u=\eta U$, $\alpha=\eta^{p-L+1}A$ transforms this expansion in an asymptotic expansion of the regularly perturbed differential equation associated to (\ref{eq-gen}), on a domain which is growing as $\eta$ tends to $0$.\\
In particular, the change of variable $\tau=\frac{t}{\eta}$ brings a local study of the canard solution in a small neighborhood of $t=0$.\vspace*{3mm}\\
\underline{Remark 2:}\\
This theorem shows that the lack of uniqueness mentionned previously is only due to the lack of control of the different interactions between the intermediary functions.\vspace*{3mm}\\
\underline{Remark 3:}\\
The asymptotic expansion presented in the last theorem are currently studied, in the complex case, by A. Fruchard and R. Sch\"afke in a recent preprint.

\section*{Acknowledgments}
I am grateful to professor \'Eric Beno\^it for the scientific framing of this work, which consists of the last part of my PhD thesis, and to professor Guy Wallet for helpful discussions as co-advisor. Moreover, I thank professor Reinhard Sch\"afke for helpful comments and improvments. This work was supported by the r\'egion Poitou-Charentes under the grant no 03/RPC-R-148.

\bibliographystyle{plain}
\bibliography{biblio-papier2}

\end{document}